\documentclass[10pt]{amsart}
\usepackage{amsthm,amsmath,amssymb,amsfonts}
\usepackage[american]{babel}
\usepackage{cite}
\usepackage{url}

\newtheorem{theorem}{Theorem}

\theoremstyle{definition}

\theoremstyle{remark}
\newtheorem{remark}[theorem]{Remark}

\begin{document}


\newcommand{\abs}[1]{\lvert#1\rvert}
\newcommand{\nt}{\trianglelefteq}
\newcommand{\NU}[1]{\mbox{\rm V}(#1)}
\newcommand{\cen}[2]{\mbox{\rm C}_{#1}(#2)}
\newcommand{\zen}[1]{\mbox{\rm Z}(#1)}
\newcommand{\nor}[2]{\mbox{\rm N}_{#1}(#2)}
\newcommand{\lara}[1]{\langle{#1}\rangle}
\newcommand{\ZZ}{\mbox{$\mathbb{Z}$}}
\newcommand{\QQ}{\mbox{$\mathbb{Q}$}}
\newcommand{\CC}{\mbox{$\mathbb{C}$}}
\newcommand{\FF}{\mbox{$\mathbb{F}$}}
\newcommand{\sQQ}{\mbox{$\scriptstyle\mathbb{Q}$}}
\newcommand{\paug}[2]{\varepsilon_{#1}(#2)}
\newcommand{\eaug}[1]{\varepsilon_{{\tt #1}}}

\newenvironment{mydescription}[1]
{\begin{list}{}
{\renewcommand{\makelabel}[1]{#1}
\settowidth{\labelwidth}{#1}
\settowidth{\leftmargin}{0pt}
\addtolength{\leftmargin}{\labelwidth}
\addtolength{\leftmargin}{-2\labelsep}
\addtolength{\leftmargin}{1pt}
}}{\end{list}}


\title[Zassenhaus conjecture for $A_{6}$]
{Zassenhaus conjecture for $\boldsymbol{A}_{\boldsymbol{6}}$}

\author{Martin Hertweck}
\address{Universit\"at Stuttgart, Fachbereich Mathematik, IGT,    
Pfaffenwald\-ring 57, 70550 Stuttgart, Germany}
\email{hertweck@mathematik.uni-stuttgart.de}

\dedicatory{Dedicated to the memory of I.~S.~Luthar$^{\dagger}$}
\thanks{$^{\dagger}$Professor Luthar died at the age of 74 
in December 2006}


\subjclass[2000]{Primary 16S34, 16U60; Secondary 20C05. \\
\hspace*{\parindent}%
{\em Subject classification.} RT}


\date{\today}

\begin{abstract}
For the alternating group $A_{6}$ of degree $6$, Zassenhaus' 
conjecture about rational conjugacy of torsion units 
in integral group rings is confirmed.
\end{abstract}

\maketitle

\section{Introduction}\label{intro}

It has proven exceedingly difficult to achieve progress in 
a number of areas revolving around units in integral group rings 
of finite groups, and this applies in particular
to the conjecture of the title:
\begin{mydescription}{(ZC1)}
\item[(ZC1)]
For a finite group $G$, every torsion unit in its integral group
ring $\ZZ G$ is conjugate to an element of $\pm G$ in the units 
of the rational group algebra $\QQ G$.
\end{mydescription}

It remains not only unsolved but also lacking in 
plausible means of either proof or counter-example.
Only a few non-solvable groups $G$ are
known to satisfy (ZC1); in this note, the conjecture is verified 
for the alternating group $A_{6}$. In 1988,
Luthar and Passi \cite{LuPa:89} started investigations in this field
when dealing with $A_{5}$, thereby introducing what is now called the
``Luthar--Passi method.'' 
A little later, Luthar and Trama \cite{LuTr:91} established
(ZC1) for the symmetric group $S_{5}$, with a large part of 
the calculations being done inside the integral group ring $\ZZ S_{5}$,
or rather a Wedderburn embedding, itself. 
Mentioning further that the covering group of $A_{5}$ can be dealt with
some additional arguments \cite{DoJuPM:97}, that accurately
describes the state of knowledge at the turn of the century.
The limitations of the original version of the Luthar--Passi method
were shown in \cite{BoHoKi:04}.
Recently, a modular version,
especially suited for application to non-solvable groups,
was given in \cite{He:05b}, and used to verify (ZC1) for the groups 
$\text{\rm PSL}(2,p)$, $p=7,11,13$. (At that point, 
the reader may wish to recall that $A_{6}\cong\text{\rm PSL}(2,9)$.)
A brief review of the method is given in Section~\ref{Sec:LP-method}.
It was also used in \cite{BoHe:06} to verify (ZC1) for 
a covering group of $S_{5}$ and for $\text{\rm GL}(2,5)$.
The method makes use of the (ordinary) character table and of modular
character tables, in a process suited for being done on a computer, 
the result being that rational conjugacy of torsion units of a given
order to group elements is either proven or not, and if not, 
at least some information about partial augmentations is obtained
(see \cite{BoKo:06,BoKo:06b,BoJeKo:06,BoKoSi:06}).

It should be remarked that from the very beginning, 
investigations on (ZC1) concentrated on solvable groups.
Here we only mention that Luthar, in collaboration with others,
initiated the study of metacyclic groups
\cite{LuBh:83,BhLu:84,LuTr:90,LuSe:98}.
For the whole class of metacyclic groups, (ZC1) has been confirmed
only very recently in \cite{He:06a}.
Other, more recent, results can be found in \cite{He:05a}
and \cite{HoKi:06}.

We briefly recall some of the necessary background.
Let $G$ be a finite group. 
For a group ring element 
$u=\sum_{g\in G}a_{g}g$ (all $a_{g}$ in $\ZZ$), its partial augmentation
$\paug{x}{u}$ with respect to an element $x$ of $G$, or rather its 
conjugacy class $x^{G}$ in $G$, is the sum $\sum_{g\in x^{G}}a_{g}$.
The augmentation of $u$ is the sum $\sum_{g\in G}a_{g}$ of all of its 
partial augmentations. Of course, we need to consider only units of
augmentation one in $\ZZ G$, which form a group we denote by $\NU{\ZZ G}$.

A criterion for a torsion unit $u$ in $\NU{\ZZ G}$ to be 
conjugate to an element of $G$ in the units of $\QQ G$, which is 
especially suited for computational purposes, is that all but one of 
the partial augmentations of every power of $u$ vanish (see 
\cite[Theorem~2.5]{MaRiSeWe:87}). 
The partial augmentations of $u$ one has to turn attention to are
limited by the following remark.

\begin{remark}\label{rem2}
Let $u$ be a torsion unit in $\NU{\ZZ G}$. Then $g\in G$ and 
$\paug{g}{u}\neq 0$ implies that the order of $g$ divides the order of $u$.
Indeed, it is well-known that then prime divisors of the order 
of $g$ divide the order of $u$ (see \cite[Theorem~2.7]{MaRiSeWe:87},
as well as \cite[Lemma~2.8]{He:05a} for an alternative
proof). Further, it was observed in \cite[Proposition~2.2]{He:05b} that
the orders of the $p$-parts of $g$ cannot exceed those of $u$.
\end{remark}

It is known (see \cite{CoLi:65}) that the order of a torsion unit
$u$ in $\NU{\ZZ G}$ is a divisor of the exponent of $G$. However, it
is not known whether some element of $G$ has the 
same order as $u$. For example, the more interesting part we will be 
concerned with when proving (ZC1) for $A_{6}$ is about the
non-existence of torsion units of order $6$ in $\NU{\ZZ A_{6}}$,
despite the fact that there are no elements of order $6$ in $A_{6}$
(see Section~\ref{Sec:u6}). 
Using the Luthar--Passi method alone is not sufficient.
In this situation, we are lucky in that
only in a rough sense `ties' between two integral representations
(affording the $5$-dimensional irreducible characters) need to be
considered.

In Section~\ref{Sec:Appl-LP}, we apply the Luthar--Passi method
to torsion units in $\NU{\ZZ A_{6}}$, with the result that all
torsion units not of order $6$ are conjugate to elements of $A_{6}$ 
in the units of $\QQ A_{6}$, and if there should exist a unit of 
order $6$, its partial augmentations are essentially unique and 
can be specified. That is a matter of routine, and was already done 
by Salim \cite{Sal:06} who calculated the coefficients of 
the linear inequalities arising from equation \eqref{eqLP2} below,
and then their solutions. We shall proceed in a straightforward way,
making only use of inequalities arising from equation \eqref{eqLP1}.

\section{The Luthar--Passi method}\label{Sec:LP-method}

We briefly survey, in an informal way, the Luthar--Passi method.
Still, let $u$ be a torsion unit in $\NU{\ZZ G}$, of order $n$ (say).
Then, conjecture (ZC1) states that there is $g\in G$ such that
$\chi(u^{i})=\chi(g^{i})$ for all irreducible characters $\chi$ of
$G$ and integers $i$. This is because complex representations are
determined (up to equivalence) by their characters, and conjugacy in
the units of $\CC G$ can be shown to be equivalent to conjugacy in
the units of $\QQ G$ (``rational conjugacy'').
It is perhaps even more natural to take the eigenvalues of representing
matrices into consideration: If an irreducible character $\chi$ is
afforded by a complex representation $D$ of $G$, write 
$\mu(\xi,u,\chi)$ for the multiplicity of an $n$-th root of unity 
$\xi$ as an eigenvalue of the matrix $D(u)$. Then $u$ is rationally
conjugate to the group element $g$ if and only if  
$\mu(\xi,u,\chi)=\mu(\xi,g,\chi)$ for all $\chi$ and $\xi$.

The same applies, mutatis mutandis, to Brauer characters.
Suppose that $u$ is a $p$-regular torsion unit for some prime $p$ 
dividing the order of $G$ (that is, $p$ does not divide the order $n$
of $u$). A Brauer character 
$\varphi$ (relative to $p$) is a complex-valued function 
associated with a modular representation of $G$ (in characteristic 
$p$), defined on the set of $p$-regular elements of $G$. 
The way Brauer defined these functions actually shows that the 
domain of $\varphi$ can be extended to the set of 
$p$-regular torsion units of $\NU{\ZZ G}$, and the decomposition
matrix gives the relation between the values of irreducible 
Brauer characters and ordinary characters on these units.
Multiplicities $\mu(\xi,u,\varphi)$ are defined in the obvious way.
Then one can formulate, in complete analogy, that
$u$ is rationally conjugate to the group element $g$ if and only if
$\varphi(u^{i})=\varphi(g^{i})$ for all irreducible Brauer characters
$\varphi$ and integers $i$, or $\mu(\xi,u,\varphi)=\mu(\xi,g,\varphi)$
for all $\varphi$ and $\xi$.

Let $\psi$ be either an irreducible ordinary character, or an 
irreducible Brauer character of $G$ relative to $p$ when $u$ is 
$p$-regular.
The multiplicities $\mu(\xi,u,\psi)$ can be computed from the values
$\psi(u^{i})$ by Fourier inversion, and taking the Galois action
into account, the formula reads---with $\zeta$ denoting a
primitive $n$-th root of unity:
\[ \mu(\xi,u,\psi) = \frac{1}{n}\sum_{d\mid n}
\text{\rm Tr}_{\sQQ(\zeta^{d})/\sQQ}(\psi(u^{d})\xi^{-d}). \]
Partial augmentations come into play by means of 
\begin{equation}\label{eqLP1}
\psi(u)=\sum_{x^{G}}\paug{x}{u}\psi(x), 
\end{equation}
where it is understood that a Brauer character
vanishes off the $p$-regular elements.
This is obvious if $\psi$ is an ordinary character, and can be seen as 
a consequence of the non-singularity of the Cartan matrix if 
$\psi$ is a Brauer character, remembering
the fact that the partial augmentations of the $p$-regular
torsion unit $u$ at $p$-singular group elements are zero.
Combining both equations, one obtains
\begin{equation}\label{eqLP2}
\mu(\xi,u,\psi) = \frac{1}{n}\sum_{x^{G}}\sum_{d\mid n}
\paug{x}{u^{d}}\text{\rm Tr}_{\sQQ(\zeta^{d})/\sQQ}(\psi(x)\xi^{-d}). 
\end{equation}
This should be seen as a linear system of equations with rational
coefficients in the unknown partial augmentations. Note that the 
$\mu(\xi,u,\psi)$ are non-negative integers, bounded above by $\psi(1)$.

\section{Application of the Luthar--Passi method}\label{Sec:Appl-LP}

Let $u$ be a torsion unit in $\NU{\ZZ A_{6}}$.
The order of $u$ is a divisor
of the exponent of $A_{6}$, which is $2^{2}\cdot 3\cdot 5$. The 
character table of $A_{6}$ is shown in Table~\ref{T1} in the form 
obtained by requiring 
{\tt CharacterTable("A6")} in GAP \cite{GAP4} (dots indicate zeros).
\begin{table}[h] 
\[
\begin{array}{c}
\begin{array}{r|rrrrrrr} \hline
& {\tt 1a} & {\tt 2a} & {\tt 3a} & {\tt 3b} & {\tt 4a} & {\tt 5a}
& {\tt 5b} \rule[-7pt]{0pt}{20pt} \\ \hline
\chi_{1a} & 1& 1& 1& 1& 1& 1& 1
\rule[0pt]{0pt}{13pt} \\
\chi_{5a} & 5& 1&  2& -1& -1& .& . \\
\chi_{5b} & 5& 1&  -1& 2& -1& .& . \\
\chi_{8a} & 8& .&  -1& -1& .& \alpha_{1} & \alpha_{2} \\
\chi_{8b} & 8& .&  -1& -1& .& \alpha_{2} & \alpha_{1} \\
\chi_{9a} & 9& 1& .& .& 1& -1& -1 \\
\chi_{10a} & 10& -2& 1& 1& .& . & . 
\rule[-7pt]{0pt}{5pt} \\ \hline
\end{array} \\
\begin{array}{rl}
\text{Irrational entries:} &
\alpha_{1}=(1+\sqrt{5})/2=
1+\nu+\nu^{4} \text{ and} \\
& \alpha_{2}=(1-\sqrt{5})/2=
1+\nu^{2}+\nu^{3}, \text{ where } 
\nu=\exp(2\pi i/5).
\rule[6pt]{0pt}{5pt} \\
\end{array}
\rule[0pt]{0pt}{25pt}
\end{array} \] 
\rule[0pt]{0pt}{5pt}
\caption{Character table of $A_{6}$}\label{T1}
\end{table}
So $k${\tt x} is the ``x''-th class of elements of order $k$, and 
$\eaug{\text{$k$}x}$ will denote the partial augmentation of $u$ with 
respect to this class. 
We will also make use of the Brauer character tables of $A_{6}$ for
the primes $2$ and $3$, shown in Table~\ref{T2} 
(which can be accessed from the GAP library through the {\tt mod}
command if the character table is required in the above form). 
\begin{table}[h] 
\[ 
\begin{array}[b]{r|r@{\hspace*{5pt}}r@{\hspace*{5pt}}r@{\hspace*{5pt}}
r@{\hspace*{5pt}}r@{\hspace*{5pt}}r} \hline
& {\tt 1a} & {\tt 3a} & {\tt 3b} & {\tt 5a} & {\tt 5b}
\rule[-7pt]{0pt}{20pt} \\ \hline
\varphi_{2,1a} & 1 &  1 &  1 &  1 &  1 \rule[0pt]{0pt}{13pt} \\
\varphi_{2,4a} & 4 &  1 & -2 & -1 & -1 \\
\varphi_{2,4b} & 4 & -2 &  1 & -1 & -1 \\
\varphi_{2,8a} & 8 & -1 & -1 & \alpha_{1} & \alpha_{2} \\
\varphi_{2,8b} & 8 & -1 & -1 & \alpha_{2} & \alpha_{1} 
\rule[-7pt]{0pt}{5pt} \\ \hline
\end{array} 
\qquad
\begin{array}[b]{r|r@{\hspace*{5pt}}r@{\hspace*{5pt}}r@{\hspace*{5pt}}
r@{\hspace*{5pt}}r@{\hspace*{5pt}}r} \hline
& {\tt 1a} & {\tt 2a} & {\tt 4a} & {\tt 5a} & {\tt 5b}
\rule[-7pt]{0pt}{20pt} \\ \hline
\varphi_{3,1a} & 1 &  1 &  1 &  1 &  1 \rule[0pt]{0pt}{13pt} \\
\varphi_{3,3a} & 3 & -1 &  1 & \alpha_{1} & \alpha_{2} \\
\varphi_{3,3b} & 3 & -1 &  1 & \alpha_{2} & \alpha_{1} \\
\varphi_{3,4a} & 4 &  . & -2 & -1 & -1 \\
\varphi_{3,9b} & 9 &  1 &  1 & -1 & -1 \rule[-7pt]{0pt}{5pt} \\
\hline
\end{array} 
\] 
\rule[0pt]{0pt}{5pt}
\vspace*{-6pt}
\caption{Brauer character tables of $A_{6}$ mod $2$ and mod $3$}
\label{T2}
\vspace*{-10pt}
\end{table}
We assume that $u\neq 1$, so $\eaug{1a}=0$ by the familiar 
Berman--Higman result. We shall show that all but one of the partial
augmentations of $u$ vanish if its order is not $6$. When $u$ is assumed
to have order $6$, the Luthar--Passi 
method is not strong enough to exclude the possibility of having
$(\eaug{2a},\eaug{3a},\eaug{3b})$ equal to $(-2,1,2)$ or $(-2,2,1)$.
Note that we have to consider the cases 
when $u$ has order $2$, $4$, $3$, $5$ (then vanishing of all but one
of the partial augmentations will be shown), when $u$ has order $10$
or $15$ (a contradiction must be reached), and when $u$ has order $6$ 
(with the final contradiction to be obtained in the next section). 
Remark~\ref{rem2} on the vanishing of certain  partial augmentations
will be used without further mention. Finally, we remark that $A_{6}$ 
has automorphisms which permute either the classes {\tt 3a} and {\tt 3b}
or the classes {\tt 5a} and {\tt 5b} (leaving the remaining classes 
fixed), which helps from repeating steps.

A subscript $\ast$ will always stand as a placeholder for all possible
subscripts.
We denote by $\Theta_{\ast}$ a representation of $G$, over a sufficiently
large field, which affords the Brauer character $\varphi_{\ast}$, and
write $\Theta_{\ast}(u)\sim\text{\rm diag}(\lambda_{1},\dotsc,
\lambda_{\varphi_{\ast}(1)})$ to indicate that $\Theta_{\ast}(u)$ 
has eigenvalues $\lambda_{1},\dotsc,\lambda_{\varphi_{\ast}(1)}$ 
(multiplicities taken into account), which will, in the sense of Brauer,
be identified with complex roots of unity.

\subsection*{When $\boldsymbol{u}$ has order $\boldsymbol{2}$.}
There is only one class of elements of order $2$,
so there is nothing to do.

\subsection*{When $\boldsymbol{u}$ has order $\boldsymbol{4}$.}
Taking augmentation gives $\eaug{2a}+\eaug{4a}=1$, so 
$\abs{-2\eaug{2a}+1}=\abs{-\eaug{2a}+\eaug{4a}}=\abs{\varphi_{3,3a}(u)}<3$
(strict inequality as 
$\abs{\varphi_{3,3a}(u^{2})}=\abs{\varphi_{3,3a}({\tt 2a})}\neq 3$).
Thus $\eaug{2a}\in\{0,1\}$, and one of $\eaug{2a}$ and $\eaug{4a}$
is zero (of course, this must be $\eaug{2a}$).

\subsection*{When $\boldsymbol{u}$ has order $\boldsymbol{3}$.}
Taking augmentation gives $\eaug{3a}+\eaug{3b}=1$. Thus 
$\abs{3\eaug{3a}-2}=\abs{\eaug{3a}-2\eaug{3b}}=
\abs{\varphi_{2,4a}(u)}\leq 4$ and
$\abs{-3\eaug{3a}+1}=\abs{-2\eaug{3a}+\eaug{3b}}=
\abs{\varphi_{2,4b}(u)}\leq 4$. So
$\eaug{3a}\in\{0,1,2\}\cap\{-1,0,1\}=\{0,1\}$, and we are done.

\subsection*{When $\boldsymbol{u}$ has order $\boldsymbol{5}$.}
Taking augmentation gives $\eaug{5a}+\eaug{5b}=1$, so 
$\varphi_{3,3a}(u)=
(1+\nu+\nu^{4})\eaug{5a}+(1+\nu^{2}+\nu^{3})\eaug{5b}=
1+(1-\eaug{5a})(\nu^{2}+\nu^{3})+\eaug{5a}(\nu+\nu^{4})$.
Since $\varphi_{3,3a}(u)$ is a sum of three fifth root of unity, 
it follows that $\eaug{5a}\in\{0,1\}$, as desired.

\medskip
By slight abuse of notation, we shall also write $k${\tt x} for a 
representative of the class $k${\tt x}.

\subsection*{When $\boldsymbol{u}$ has order $\boldsymbol{10}$.}
Taking augmentation gives $\eaug{2a}+\eaug{5a}+\eaug{5b}=1$. Thus
$\varphi_{3,3a}(u)=
-1+\eaug{5a}(3+\sqrt{5})/2+\eaug{5b}(3-\sqrt{5})/2$.
We can assume without lost of generality that $u^{6}$ is rationally
conjugate to {\tt 5a}. Also $u^{5}$ is rationally
conjugate to {\tt 2a}. We have
$\Theta_{3,3a}({\tt 5a})\sim\text{\rm diag}(1,\nu,\nu^{4})$ and
$\Theta_{3,3a}({\tt 2a})\sim\text{\rm diag}(1,-1,-1)$. 
The eigenvalues of $\Theta_{3,3a}(u)$ are products of the
diagonal entries of the diagonal matrices just depicted,
always taken one from each matrix,
in such a way that all entries are involved.
Since $\varphi_{3,3a}(u)$ is real, it follows that 
$\Theta_{3,3a}(u)\sim\text{\rm diag}(1,-\nu,-\nu^{4})$ and
$\varphi_{3,3a}(u)=(3-\sqrt{5})/2$. Comparing the given values for 
$\varphi_{3,3a}(u)$ results in an equation which has no integral 
solution (reduce modulo $3$), a contradiction.

\subsection*{When $\boldsymbol{u}$ has order $\boldsymbol{15}$.}
Then we can assume without lost of generality that $u^{10}$ is rationally
conjugate to {\tt 3a} and that $u^{6}$ is rationally
conjugate to {\tt 5a}. Setting $\zeta=\exp(2\pi i/3)$, we have
$\Theta_{2,4a}({\tt 3a})\sim\text{\rm diag}(1,1,\zeta,\zeta^{2})$ and
$\Theta_{2,4a}({\tt 5a})\sim\text{\rm diag}(\nu,\nu^{2},\nu^{3},\nu^{4})$.
The eigenvalues of $\Theta_{2,4a}(u)$ are the products of the diagonal
entries of these diagonal matrices, taken in a suitable order. But there 
is no possibility to choose an ordering such that the sum of the products
becomes rational, a contradiction.

\subsection*{When $\boldsymbol{u}$ has order $\boldsymbol{6}$.}
Taking augmentation gives $\eaug{2a}+\eaug{3a}+\eaug{3b}=1$.
Thus $\chi_{5a}(u)=1+\eaug{3a}-2\eaug{3b}$ and
$\chi_{5b}(u)=1-2\eaug{3a}+\eaug{3b}$. Now, $u^{3}$ is rationally
conjugate to {\tt 2a}, and we can assume without lost of generality 
that $u^{2}$, and hence also $u^{4}$, is rationally conjugate to {\tt 3b}. 
Write $D_{5a}$ and $D_{5b}$ for ordinary
representations affording $\chi_{5a}$ and $\chi_{5b}$, respectively.
We have $D_{5a}({\tt 2a})\sim\text{\rm diag}(1,1,1,-1,-1)$ and
$D_{5a}({\tt 3b})\sim\text{\rm diag}(1,\zeta,\zeta^{2},\zeta,\zeta^{2})$.
The eigenvalues of $D_{5a}(u)$ are products of the diagonal entries
of these diagonal matrices, suitably taken. Their sum $\chi_{5a}(u)$ 
is rational, which forces 
\begin{equation}\label{eq1}
D_{5a}(u)\sim\text{\rm diag}(1,\zeta,\zeta^{2},-\zeta,-\zeta^{2})
\end{equation}
and $\chi_{5a}(u)=1$.
Thus $\eaug{3a}=2\eaug{3b}$ and $\chi_{5b}(u)=1-3\eaug{3b}$.
We have $D_{5b}({\tt 2a})\sim\text{\rm diag}(1,1,1,-1,-1)$ and
$D_{5b}({\tt 3b})\sim\text{\rm diag}(1,1,1,\zeta,\zeta^{2})$. In the 
same way as before, we obtain that either
\begin{equation}\label{eq2}
D_{5b}(u)\sim\text{\rm diag}(1,-1,-1,\zeta,\zeta^{2})
\end{equation}
and $\chi_{5b}(u)=-2$, or $\chi_{5b}(u)=4$.
If $\chi_{5b}(u)=4$, then $\eaug{3b}=-1$, $\eaug{3a}=-2$ and  
$\eaug{2a}=4$; but then $\chi_{10a}(u)=-11$, a contradiction.
Hence $\chi_{5b}(u)=-2$, and $\eaug{3b}=1$, $\eaug{3a}=2$ and  
$\eaug{2a}=-2$, one of the possibilities mentioned at the beginning.

\section{Non-existence of a unit of order $6$ in $\NU{\ZZ A_{6}}$}\label{Sec:u6}

Finally, we assume that there exists a unit $u$ of order $6$ in 
$\NU{\ZZ A_{6}}$, and are looking for a contradiction. 
According to the previous section, we can assume that $u$ 
has partial augmentations $(\eaug{2a},\eaug{3a},\eaug{3b})$ equal to 
$(-2,2,1)$, that $u^{3}$ is rationally conjugate to {\tt 2a} and that 
$u^{2}$ is rationally conjugate to {\tt 3b}. In particular, the
character values of all powers of $u$ are known,
from which the eigenvalues of representing matrices can be computed.

We shall take advantage only of representations $D_{5a}$ and $D_{5b}$
affording the irreducible characters $\chi_{5a}$ and $\chi_{5b}$, for 
which the eigenvalues of $D_{5a}(u)$ and $D_{5b}(u)$ are already 
given by \eqref{eq1} and \eqref{eq2}, taking into account that
$\chi_{5a}$ and $\chi_{5b}$ have the same
$3$-modular constituents, namely the principal character and the
$3$-modular character of degree $4$ (as can be seen from 
Tables~\ref{T1} and \ref{T2}).

Set $R=\ZZ_{(3)}$, the localization of $\ZZ$ at the prime $3$,
and $k=R/3R$ (the prime field of characteristic $3$). 
For an $RA_{6}$-lattice $L$, we let $\bar{L}$ stand for the 
$kA_{6}$-module $L/3L$, and we also let $x\mapsto\bar{x}$ 
denote the natural homomorphism from $RA_{6}$ to $kA_{6}$.

There are two doubly transitive actions of $A_{6}$ on a set
with six elements, corresponding to the two conjugacy classes
of subgroups of $A_{6}$ which are isomorphic to $A_{5}$, and 
$D_{5a}$ and $D_{5b}$ can be chosen to be the associated
deleted permutation representations.  
Denote by $L_{a}$ and $L_{b}$ the non-trivial factor lattices
of the $RA_{6}$-permutation lattices arising from these actions,
so that a matrix representation $D_{5\ast}$ affording $\chi_{5\ast}$
can be obtained with respect to an $R$-basis of $L_{\ast}$.

Set $e=(1+u^{3})/2$ and $f=1-e$, so $1=e+f$ is an idempotent
decomposition in $RA_{6}$. Let $w=e+fu^{2}\in RA_{6}$. 
From \eqref{eq2} one sees that $D_{5b}(w)$ is the identity matrix.
In particular, $\bar{w}$ acts as the identity on the composition
factors of $\bar{L}_{b}$, which are also the composition
factors of $\bar{L}_{a}$, as already noted.

Now we turn our attention to the other representation.
Set $L=L_{a}$. Note that the simple $kA_{6}$-module $\tilde{M}$
of dimension $4$ is a submodule of $\bar{L}$, namely 
$\tilde{M}=\bar{L}\text{\rm I}(kA_{6})$ where $\text{\rm I}(kA_{6})$
denotes the augmentation ideal of $kA_{6}$.
(If we set $M=L\text{\rm I}(RA_{6})$, then $M$ is a sublattice
of $L$ of index $3$ and $\tilde{M}$ is the image of $M$ under the
map $L\rightarrow\bar{L}$.) 
We have $L=Le\oplus Lf$. From \eqref{eq1} one sees that 
$Lf$ is an $R$-lattice of rank $2$ on which $w$ acts 
nontrivially by multiplication with $u^{2}$, a unit of order $3$.
But $\overline{\!Lf}=\bar{L}\bar{f}\subseteq\tilde{M}$, so
$\bar{w}$ acts trivially on $\overline{\!Lf}$, by the 
preceding paragraph. So we have obtained a contradiction, as
the kernel of the natural map 
$\text{\rm GL}_{2}(R)\rightarrow\text{\rm GL}_{2}(k)$
has no $3$-torsion. (A particular case of a well-known lemma
due to Minkowski; see, for example, the beautiful discussion in  
\cite[5.2]{GuLo:06}.)



\providecommand{\bysame}{\leavevmode\hbox to3em{\hrulefill}\thinspace}
\providecommand{\href}[2]{#2}

\end{document}